\documentclass[a4paper,12pt,twoside,leqno,final]{amsart}
\usepackage{amsmath}
\usepackage{amssymb}

\newtheorem{thm}{Theorem}[section]

\newtheorem{cor}[thm]{Corollary}

\newcommand{\C}{{\mathbb C}}
\newcommand{\D}{{\mathbb D}}

\newcommand{\T}{{\mathbb T}}

\newcommand{\f}{\frac}
\newcommand{\ov}{\overline}
\newcommand{\al}{\alpha}

\newcommand{\la}{\lambda}
\newcommand{\ze}{\zeta}
\renewcommand{\th}{\theta}

\newcommand{\ph}{\varphi}
\newcommand{\om}{\omega}

\numberwithin{equation}{section}

\title[A reverse Schwarz--Pick inequality]
{A reverse Schwarz--Pick inequality}
\author{Konstantin M. Dyakonov}
\address{ICREA and Universitat de Barcelona, Departament de Matem\`atica 
Aplicada i An\`alisi, Gran Via 585, E-08007 Barcelona, Spain}
\email{konstantin.dyakonov@icrea.cat}
\thanks{Supported in part by grant MTM2011-27932-C02-01 from El Ministerio de Ciencia 
e Innovaci\'on (Spain) and grant 2009-SGR-1303 from AGAUR (Generalitat de Catalunya).}

\begin{document}
\begin{abstract}
We prove a kind of \lq\lq reverse Schwarz--Pick lemma" for holomorphic self-maps of the disk. The result 
becomes especially clear-cut for inner functions and casts new light on their derivatives. 
\end{abstract}

\maketitle

\section{Introduction and main result} 

We write $\D$ for the disk $\{z\in\C:|z|<1\}$, $\T$ for its boundary, and $m$ for the normalized arclength 
measure on $\T$; thus $dm(\ze)=(2\pi)^{-1}|d\ze|$. Further, $H^\infty$ will denote the algebra of bounded 
holomorphic functions on $\D$, equipped with the usual supremum norm $\|\cdot\|_\infty$. 

\par The classical Schwarz lemma and its invariant version, known as the Schwarz--Pick lemma 
(see \cite[Chapter I]{G}), lie at the heart of function theory on the disk. The Schwarz--Pick lemma 
tells us, in particular, that every function $\ph\in H^\infty$ with $\|\ph\|_\infty\le1$ satisfies 
\begin{equation}\label{eqn:schpick}
|\ph'(z)|\le\f{1-|\ph(z)|^2}{1-|z|^2},\qquad z\in\D.
\end{equation}
Moreover, equality holds at some---or each---point of $\D$ if and only if $\ph$ is either 
a {\it M\"obius transformation} (i.e., has the form 
$$z\mapsto\la\f{z-a}{1-\ov az}$$
for some $\la\in\T$ and $a\in\D$) or a constant of modulus 1. 
It seems natural to ask how large the gap between the two quantities 
in \eqref{eqn:schpick} may be in general. To be more precise, we would like to know 
whether some sort of reverse estimate has a chance to hold, i.e., whether the quotient 
$$\f{1-|\ph(z)|^2}{1-|z|^2}=:Q_\ph(z)$$
admits an upper bound in terms of $\ph'$. And if it does, what is the appropriate $\ph'$-related 
majorant for that quantity? 

\par Of course, there are limits to what can be expected. Assume, from now on, 
that $\ph$ is a {\it nonconstant} $H^\infty$-function of norm at most $1$. 
First of all, while $\ph'$ may happen to vanish or become arbitrarily small at certain points 
of $\D$, it follows from Schwarz's lemma that $Q_\ph$ is always bounded away from zero; in fact, 
\begin{equation}\label{eqn:estfrombel}
Q_\ph(z)\ge\f{1-|\ph(0)|}{1+|\ph(0)|},\qquad z\in\D.
\end{equation}
As another example, consider the case where $\ph$ is analytic in a neighborhood 
of some closed arc $\gamma\subset\T$ and satisfies $\sup_{\ze\in\gamma}|\ph(\ze)|<1$. 
In this case, too, the two sides of \eqref{eqn:schpick} exhibit different types of behavior 
as $z$ approaches $\gamma$. Indeed, $|\ph'(z)|$ is then well-behaved, whereas $Q_\ph(z)$ blows up 
like a constant times $(1-|z|)^{-1}$. Thus, the ratio $Q_\ph(z)/|\ph'(z)|$ is sometimes huge. 

\par On the other hand, suppose that $\ph$ has an {\it angular derivative} (in the sense of 
Carath\'eodory) at a point $\ze\in\T$. This means, by definition, that $\ph$ and $\ph'$ both 
have nontangential limits at $\ze$ and, once we agree to denote the two limits by $\ph(\ze)$ 
and $\ph'(\ze)$, the former of these satisfies $|\ph(\ze)|=1$. The classical Julia--Carath\'eodory 
theorem (see \cite[Chapter VI]{B}, \cite[Chapter I]{Car} or \cite[Chapter VI]{Sar2}) 
asserts that this happens if and only if 
$$\liminf_{z\to\ze}\f{1-|\ph(z)|}{1-|z|}<\infty,$$ 
a condition that can be further rewritten as 
\begin{equation}\label{eqn:liminf}
\liminf_{z\to\ze}Q_\ph(z)<\infty.
\end{equation}
And if any of these holds, the theorem tells us also that $\ph'(\ze)$ coincides with the limit of the 
difference quotient 
$$\f{\ph(z)-\ph(\ze)}{z-\ze}$$ 
as $z\to\ze$ nontangentially, while the nontangential limit of $Q_\ph(z)$ equals $|\ph'(\ze)|$. In 
addition, this last number agrees with the value of the (unrestricted) $\liminf$ in \eqref{eqn:liminf}. 

\par Consequently, if $\ph$ happens to possess an angular derivative $\ph'(\ze)$ at every point $\ze$ of 
a set $\mathcal E\subset\T$, then the two sides of \eqref{eqn:schpick} have the same boundary values 
on $\mathcal E$. We may therefore expect the two quantities to be reasonably close near $\mathcal E$, 
so a certain \lq\lq reverse Schwarz--Pick type estimate" is likely to hold on a suitable region 
of $\D$ adjacent to $\mathcal E$. And the more massive $\mathcal E$ is, the stronger should our 
reverse inequality become. Our main result, Theorem \ref{thm:mainres} below, provides such an estimate 
under the hypotheses that the \lq\lq good" set $\mathcal E$ is (Lebesgue) measurable and 
the function $\log|\ph'|$ is integrable on $\mathcal E$; of course, only the case $m(\mathcal E)>0$ 
is of interest. 

\par Before stating the result, let us recall that the {\it harmonic measure} $\om_z$ associated with 
a point $z\in\D$ is given by 
$$d\om_z(\ze)=\f{1-|z|^2}{|\ze-z|^2}\,dm(\ze),\qquad\ze\in\T.$$ 
For a measurable set $E\subset\T$, the quantity $\om_z(E)=\int_Ed\om_z$ is thus the value at $z$ 
of the harmonic extension (into $\D$) of the characteristic function $\chi_E$; this quantity can 
be roughly thought of as the normalized angle at which $E$ is seen from $z$. We also recall that 
if $h$ is a nonnegative function on $\T$ with $\log h\in L^1(\T)$, then 
$$\mathcal O_h(z):=\exp\left\{\int_\T\f{\ze+z}{\ze-z}\log h(\ze)\,dm(\ze)\right\},
\qquad z\in\D,$$
is a holomorphic function on $\D$ whose modulus has nontangential boundary values $h$ almost everywhere 
on $\T$. In fact, 
$$|\mathcal O_h(z)|=\exp\left\{\int_\T\log h\,d\om_z\right\},$$ 
whence the preceding statement follows. This function $\mathcal O_h$ is known as the {\it outer function 
with modulus $h$}. Functions of the form $\la\mathcal O_h$, with $\la\in\T$ and $h$ as above, will also be 
referred to as outer; see \cite[Chapter II]{G} for a further discussion of outer functions and their properties. 

\begin{thm}\label{thm:mainres} Let $\ph\in H^\infty$ be a nonconstant function with $\|\ph\|_\infty=1$, 
and let $\mathcal E$ be a measurable subset of $\T$ such that $\ph$ has an angular derivative 
almost everywhere on $\mathcal E$. Assume that $\log|\ph'|\in L^1(\mathcal E)$ and write $\mathcal O$ 
for the outer function with modulus $|\ph'|\chi_\mathcal E+\chi_{\T\setminus\mathcal E}$. Then 
\begin{equation}\label{eqn:mainest}
Q_\ph(z)\le|\mathcal O(z)|\left\{\f1{1-\om_z(\mathcal E)}\cdot\f{1+|z|}{1-|z|}\right\}^{1-\om_z(\mathcal E)}
\end{equation}
for all $z\in\D$. 
\end{thm}

\par The first factor on the right being 
$$|\mathcal O(z)|=\exp\left\{\int_\mathcal E\log|\ph'|\,d\om_z\right\},$$ 
we may indeed view \eqref{eqn:mainest} as a reverse Schwarz--Pick inequality, since it provides 
us with an upper bound for $Q_\ph$ in terms of $|\ph'|$. If, in addition, $|\ph'|$ happens to be 
essentially bounded on $\mathcal E$, then we clearly have 
\begin{equation}\label{eqn:triv}
|\mathcal O(z)|\le\left\|\ph'\right\|_{\infty,\mathcal E}^{\om_z(\mathcal E)},
\end{equation}
where $\|\cdot\|_{\infty,\mathcal E}=\|\cdot\|_{L^\infty(\mathcal E)}$. Combining \eqref{eqn:mainest} with 
\eqref{eqn:triv} and with the elementary fact that 
$$\sup\left\{t^{-t}:\,0<t<1\right\}=e^{1/e}$$ 
leads to the weaker, but perhaps simpler, estimate 
\begin{equation}\label{eqn:cormainest}
Q_\ph(z)\le e^{1/e}\left\|\ph'\right\|_{\infty,\mathcal E}^{\om_z(\mathcal E)}
\left(\f{1+|z|}{1-|z|}\right)^{1-\om_z(\mathcal E)},\qquad z\in\D.
\end{equation}

\par Both \eqref{eqn:mainest} and \eqref{eqn:cormainest} reflect the influence of the \lq\lq good" set 
$\mathcal E$ and of the \lq\lq bad" set $\T\setminus\mathcal E$, depending on the location of $z$, in 
the spirit of Nevanlinna's {\it Zweikonstantensatz} (the two constants theorem). The latter 
supplies a sharp bound on $|f(z)|$ for an $H^\infty$-function $f$ whose modulus is bounded by two given 
constants on two mutually complementary subsets of the boundary; see, e.\,g., \cite[Chapter VIII]{E}. 

\par The next section contains some applications of Theorem \ref{thm:mainres} to inner functions, while 
the proof of the theorem is given in Section 3.

\section{Inner functions and their derivatives}

Recall that a function $\th\in H^\infty$ is said to be {\it inner} if $\lim_{r\to1^-}|\th(r\ze)|=1$ for 
almost all $\ze\in\T$. Also involved in what follows is the {\it Nevanlinna class} $\mathcal N$, defined 
as the set of all holomorphic functions $f$ on $\D$ that satisfy 
$$\sup_{0<r<1}\int_\T\log^+|f(r\ze)|\,dm(\ze)<\infty.$$ 
Equivalently, $\mathcal N$ is formed by the ratios $u/v$ with $u,v\in H^\infty$ and with $v$ zero-free 
on $\D$; see \cite[Chapter II]{G}. Imposing the additional restriction that $v$ be outer, one arrives at 
the definition (or a characterization) of the {\it Smirnov class} $\mathcal N^+$. 
\par Now, if $\th$ is an inner function with $\th'\in\mathcal N$, then $\th$ has an angular derivative a.\,e. 
on $\T$ and $\log|\th'|\in L^1(\T)$. Therefore, when applying Theorem \ref{thm:mainres} to $\ph=\th$, we 
may take $\mathcal E=\T$. This yields the following result. 

\begin{cor}\label{cor:inn1} Let $\th$ be a nonconstant inner function with $\th'\in\mathcal N$, and let 
$\mathcal O=\mathcal O_{|\th'|}$ be the outer factor of $\th'$ (i.\,e., the outer function with 
modulus $|\th'|$ on $\T$). Then 
\begin{equation}\label{eqn:innest}
Q_\th(z)\le|\mathcal O(z)|,\qquad z\in\D.
\end{equation}
\end{cor}

A similar estimate can be found in \cite{DSpb}. We also mention that there is an alternative route to 
Corollary \ref{cor:inn1} via subharmonicity, which hinges on Lemma 1.1 from \cite{HSpb}; this approach 
was kindly brought to my attention by Haakan Hedenmalm. 

As a consequence of the preceding result, we now derive an amusing characterization of M\"obius transformations. 

\begin{cor}\label{cor:inn2} Given a nonconstant inner function $\th$ with $\th'\in\mathcal N$, the 
following are equivalent. 
\par{\rm(i)} $\th$ is a M\"obius transformation. 
\par{\rm(ii)} $\th'$ is an outer function. 
\par{\rm(iii)} There is a nondecreasing function $\eta:(0,\infty)\to(0,\infty)$ such that 
\begin{equation}\label{eqn:esteta}
\eta\left(Q_\th(z)\right)\le|\th'(z)|,\qquad z\in\D.
\end{equation}
\end{cor}

Before proving this, we recall that an inner function $\th$ with $\th'\in\mathcal N$ will automatically 
have $\th'$ in $\mathcal N^+$, a fact established by Ahern and Clark in \cite{AC}. 

\medskip\noindent{\it Proof of Corollary \ref{cor:inn2}.} The (i)$\implies$(ii) part is 
straightforward, while the converse follows from \eqref{eqn:innest}. Indeed, if $\th'$ is outer, 
then $|\mathcal O(z)|$ on the right-hand side of \eqref{eqn:innest} coincides with $|\th'(z)|$. 
Combining this with the Schwarz--Pick inequality \eqref{eqn:schpick}, where we put $\ph=\th$, gives 
\begin{equation}\label{eqn:speq}
|\th'(z)|=Q_\th(z),\qquad z\in\D.
\end{equation}
Thus, the current choice of $\ph$ ensures equality in \eqref{eqn:schpick}, so $\th$ must be a M\"obius 
transformation. 
\par Now that (i) and (ii) are known to be equivalent, it suffices to show that 
(i)$\implies$(iii)$\implies$(ii). The first of these implications is obvious, since every M\"obius 
transformation satisfies \eqref{eqn:speq}, so that \eqref{eqn:esteta} holds with $\eta(t)=t$. Finally, 
assuming (iii) and using \eqref{eqn:estfrombel} with $\ph=\th$, we deduce that $|\th'(z)|$ is bounded 
away from zero on $\D$, whence $1/\th'\in H^\infty$. Because $\th'\in\mathcal N^+$, it follows that $\th'$ is 
an outer function, and we arrive at (ii). 
\quad\qed

\medskip The (ii)$\implies$(i) part of Corollary \ref{cor:inn2} can be rephrased by saying that the 
derivative $\th'$ of a non-M\"obius inner function $\th$ is never outer, as long as it is in $\mathcal N$. 
Some special cases of this statement have been known. In particular, it was proved by Ahern and Clark 
(see \cite[Corollary 4]{AC}) that, within the current class of $\th$'s, the inner part of $\th'$ will 
be nontrivial provided that $\th$ has a singular factor (because such factors are actually inherited 
by $\th'$). On the other hand, if $\th$ is a finite Blaschke product with at least two zeros, 
then $\th'$ is known to have zeros in $\D$ (see \cite{W} for more information on the location 
of these), so it is clear, once again, that $\th'$ is non-outer. The same conclusion is obviously 
true for those inner functions $\th$ which have multiple zeros in $\D$. 
\par The interesting case is, therefore, that of an infinite Blaschke product with simple 
zeros. Note that if $B$ is such a Blaschke product, then, unlike in the rational 
case, $B'$ may well be zero-free on $\D$. For instance, this happens for 
$$B_\al(z):=\f{S(z)-\al}{1-\ov\al S(z)},$$ 
where $S$ is the \lq\lq atomic" singular inner function given by 
$$S(z):=\exp\left(\f{z+1}{z-1}\right)$$
and $\al$ is a point in $\D\setminus\{0\}$. One easily verifies that $B_\al$ is 
indeed a Blaschke product, while the inner factor of $B'_\al$ is $S$. 
\par We conclude this section with a question. Let $\mathfrak I$ stand for the set of nonconstant 
inner functions. {\it Which inner functions occur as inner factors (and/or divisors of such factors) 
for functions in $\mathcal N\cap\{\th':\th\in\mathfrak I\}$?} One immediate observation is that 
if $I$ is inner and $I'\in\mathcal N$, then $I$ divides the inner part of $(I^2)'$, and this last 
function is in $\mathcal N$.

\section{Proof of Theorem \ref{thm:mainres}} 

For all $z\in\D$ and almost all $\ze\in\mathcal E$, Julia's lemma (see \cite[p.\,41]{G}) yields 
\begin{equation}\label{eqn:julia}
\f{|\ph(\ze)-\ph(z)|^2}{1-|\ph(z)|^2}\le|\ph'(\ze)|\cdot\f{|\ze-z|^2}{1-|z|^2}, 
\end{equation}
or equivalently, 
\begin{equation}\label{eqn:juliabis}
\f{1-|z|^2}{1-|\ph(z)|^2}\cdot\left|\f{1-\ov{\ph(z)}\ph(\ze)}{1-\ov z\ze}\right|^2\le|\ph'(\ze)| 
\end{equation}
(recall that $|\ph(\ze)|=1$ whenever $\ph$ has an angular derivative at $\ze$). 
Keeping $z\in\D$ fixed for the rest of the proof, we now introduce the $H^\infty$-function 
\begin{equation}\label{eqn:fzphi}
F_z(w):=\f{1-|z|^2}{1-|\ph(z)|^2}\cdot\left(\f{1-\ov{\ph(z)}\ph(w)}{1-\ov zw}\right)^2
\end{equation}
and go on to rewrite \eqref{eqn:juliabis} in the form 
\begin{equation}\label{eqn:estone}
|F_z(\ze)|\le|\ph'(\ze)|,\qquad\ze\in\mathcal E.
\end{equation}
Next, we define $G_z$ to be the outer function with modulus 
$$\left|G_z(\ze)\right|=|\ph'(\ze)|\cdot\chi_\mathcal E(\ze)
+|F_z(\ze)|\cdot\chi_{\widetilde{\mathcal E}}(\ze),\qquad \ze\in\T,$$ 
where $\widetilde{\mathcal E}:=\T\setminus\mathcal E$, and note that 
\begin{equation}\label{eqn:estoncirc}
|F_z(\ze)|\le|G_z(\ze)|,\qquad\ze\in\T.
\end{equation} 
Indeed, for $\ze\in\mathcal E$ this last inequality coincides with \eqref{eqn:estone}, while 
for $\ze\in\widetilde{\mathcal E}$ it reduces to an obvious equality. 
\par Since $G_z$ is outer, \eqref{eqn:estoncirc} implies a similar estimate on $\D$, that is, 
$$|F_z(w)|\le|G_z(w)|,\qquad w\in\D.$$
In particular, setting $w=z$, we obtain 
\begin{equation}\label{eqn:estatz}
|F_z(z)|\le|G_z(z)|.
\end{equation}
It is clear from \eqref{eqn:fzphi} that 
\begin{equation}\label{eqn:lefthand}
|F_z(z)|=F_z(z)=\f{1-|\ph(z)|^2}{1-|z|^2}=Q_\ph(z),
\end{equation}
and we take further efforts to estimate $|G_z(z)|$. 
\par We have 
\begin{equation}\label{eqn:righthand}
\log|G_z(z)|=\int_\T\log|G_z(\ze)|\,d\om_z(\ze)=I_1(z)+I_2(z),
\end{equation}
where 
\begin{equation}\label{eqn:int1}
I_1(z):=\int_\mathcal E\log|\ph'(\ze)|\,d\om_z(\ze)=\log|\mathcal O(z)|
\end{equation}
and 
\begin{equation}\label{eqn:int2}
I_2(z):=\int_{\widetilde{\mathcal E}}\log|F_z(\ze)|\,d\om_z(\ze).
\end{equation}
The function $t\mapsto\log t$ being concave for $t>0$, we find that 
\begin{equation}\label{eqn:conc}
\begin{aligned}
I_2(z)&=\om_z(\widetilde{\mathcal E})\cdot\int_{\widetilde{\mathcal E}}\log|F_z(\ze)|\,
\f{d\om_z(\ze)}{\om_z(\widetilde{\mathcal E})}\\
&\le\om_z(\widetilde{\mathcal E})\cdot\log\left\{\f1{\om_z(\widetilde{\mathcal E})}
\int_{\widetilde{\mathcal E}}|F_z(\ze)|\,d\om_z(\ze)\right\}.
\end{aligned}
\end{equation}
We proceed by observing that 
\begin{equation}\label{eqn:fzze}
|F_z(\ze)|\le\f{1+|z|}{1-|z|}\cdot\f{|1-\ov{\ph(z)}\ph(\ze)|^2}{1-|\ph(z)|^2},\qquad\ze\in\T. 
\end{equation}
Furthermore, 
\begin{equation*}
\begin{aligned}
\int_\T|1-\ov{\ph(z)}\ph(\ze)|^2\,d\om_z(\ze)&=
\int_\T\left\{1-2\,\text{\rm Re}\left(\ov{\ph(z)}\ph(\ze)\right)+|\ph(z)|^2|\ph(\ze)|^2\right\}d\om_z(\ze)\\
&\le\int_\T\left\{1-2\,\text{\rm Re}\left(\ov{\ph(z)}\ph(\ze)\right)+|\ph(z)|^2\right\}d\om_z(\ze)\\
&=1-|\ph(z)|^2.
\end{aligned}
\end{equation*}
(Here, the last step consists in integrating a harmonic function against $d\om_z$, so the output 
is the function's value at $z$). In conjunction 
with \eqref{eqn:fzze}, this gives 
$$\int_\T|F_z(\ze)|\,d\om_z(\ze)\le\f{1+|z|}{1-|z|},$$
whence {\it a fortiori} 
\begin{equation}\label{eqn:taburetka}
\int_{\widetilde{\mathcal E}}|F_z(\ze)|\,d\om_z(\ze)\le\f{1+|z|}{1-|z|}.
\end{equation}
Plugging \eqref{eqn:taburetka} into \eqref{eqn:conc}, we now get 
\begin{equation}\label{eqn:int2fin}
I_2(z)\le\om_z(\widetilde{\mathcal E})\log\left\{\f1{\om_z(\widetilde{\mathcal E})}
\cdot\f{1+|z|}{1-|z|}\right\}.
\end{equation}
\par Finally, we combine \eqref{eqn:righthand} with \eqref{eqn:int1} and \eqref{eqn:int2fin} to infer that 
$$\log|G_z(z)|\le\log|\mathcal O(z)|+\om_z(\widetilde{\mathcal E})\log\left\{\f1{\om_z(\widetilde{\mathcal E})}
\cdot\f{1+|z|}{1-|z|}\right\}$$
and hence 
\begin{equation}\label{eqn:unitaz}
|G_z(z)|\le|\mathcal O(z)|\left\{\f1{\om_z(\widetilde{\mathcal E})}
\cdot\f{1+|z|}{1-|z|}\right\}^{\om_z(\widetilde{\mathcal E})}.
\end{equation}
This done, a juxtaposition of \eqref{eqn:estatz}, \eqref{eqn:lefthand} and \eqref{eqn:unitaz} yields 
the required estimate \eqref{eqn:mainest} and completes the proof. 
\quad\qed

\medskip\noindent{\it Remark.} The function $F_z$, as defined by \eqref{eqn:fzphi}, can be written in the form 
$$F_z(w)=\f{k^2_{\ph,z}(w)}{k_{\ph,z}(z)},$$
where 
$$k_{\ph,z}(w):=\f{1-\ov{\ph(z)}\ph(w)}{1-\ov zw}$$
is the reproducing kernel for the de Branges--Rovnyak space $\mathcal H(\ph)$; see \cite{Sar2}. Even though 
this space does not show up in our proof, the appearance of its kernel functions may not be incidental. It 
should be mentioned that Hilbert space methods have been previously employed, in the $\mathcal H(\ph)$ setting, 
in connection with generalized Schwarz--Pick inequalities \cite{AR} and with the Julia--Carath\'eodory theorem 
on angular derivatives \cite{Sar1, Sar2}.

\end{document}